\documentclass[12pt, leqno]{article}
\usepackage{amsmath}

\usepackage{oxford3}
\usepackage{amssymb, amsmath, url, amsthm, amssymb,graphicx}
\usepackage[T1]{fontenc}
\usepackage[latin1]{inputenc}

\newtheorem{@capt}{\sc Table}[section]
\newenvironment{capt}{\begin{@capt}\rm}{\end{@capt}}
\newtheorem{@capf}{\sc Figure}[section]
\newenvironment{capf}{\begin{@capf}\rm}{\end{@capf}}

\newcommand \fg{\varphi}

\newtheorem{Theorem}{Theorem}[section]

\newcommand \Gg{\Gamma}

\def\convd{\stackrel{\mbox{$\scriptstyle d$}}{\rightarrow}}

\def\convP{\stackrel{\mbox{$\scriptstyle P$}}{\rightarrow}}
\newcommand \ga{\gamma}

\renewcommand{\baselinestretch}{1.1}
\newcommand{\beq}{\begin{equation}}  
 \newcommand{\eeq}{\end{equation}}
\newcommand\bzero{\mbox{\boldmath${0}$}}

 \newcommand \la{\lambda}
\newcommand{\bd}{\begin{displaymath}}  
 \newcommand{\ed}{\end{displaymath}}
\newcommand \dg{\delta}

\newcommand\bDg{\mbox{\boldmath${\Delta}$}}

\newcommand\eg{\varepsilon}

\newcommand\lnorm{\left | \left |}
\newcommand\rnorm{\right  | \right |}
\newcommand\lp{\left (}
\newcommand\rp{\right )}
\newcommand\lb{\left [}
\newcommand\rb{\right ]}
\newcommand\lmo{\left |}
\newcommand\rmo{\right |}
\newcommand\lbr{\left \{}
\newcommand\rbr{\right \}}
\newcommand\lip{\left \langle}
\newcommand\rip{\right \rangle}

\newcommand\ba{{\bf a}}

\newcommand\bQ{{\bf Q}}

\newcommand\hla{{\hat \la}}

\newcommand\hc{{\hat c}}

\newcommand\dint{\int\!\!\!\int}

\numberwithin{equation}{section}
\numberwithin{Theorem}{section}
\numberwithin{Remark}{section}
\numberwithin{Assumption}{section}
\numberwithin{Lemma}{section}

\begin{document}

\title{Estimation of the mean of functional time series and
       a two sample problem}
\author{
Lajos Horv{\'a}th\\ \normalsize University of Utah, Salt Lake City, USA
\and
Piotr Kokoszka\thanks{ {\it Address for correspondence}:
Piotr Kokoszka, Department of Mathematics and Statistics,
Utah State University, 3900 Old Main Hill,
Logan, UT 84322-3900, USA.\ \
E-mail: Piotr.Kokoszka@usu.edu}
\\ \normalsize Utah State University, Logan, USA
\and
Ron Reeder \\ \normalsize University of Utah, Salt Lake City, USA
}
\date{}
\maketitle

\begin{abstract}
\noindent
This paper is concerned with inference based on the mean function of a
functional time series, which is defined as a collection of curves
obtained by splitting a continuous time record, e.g.  into daily or
annual curves. We develop a normal approximation for the functional
sample mean, and then focus on the estimation of the asymptotic
variance kernel. Using these results, we develop and asymptotically
justify a testing procedure for the equality of means in two
functional samples exhibiting temporal dependence.  Evaluated by means
of a simulations study and application to real data sets, this two
sample procedure enjoys good size and power in finite samples.  We
provide the details of its numerical implementation.

\vspace{8mm}


\vspace{2mm}

\noindent {\em Abbreviated Title:} Functional time series

\end{abstract}

\newpage

\section{Introduction} \label{s:I}
Functional time series form a class
of data structures which occurs in many applications, but several
important aspects of estimation and testing for such data have not
received as much attention as for functional data derived from
randomized experiments. In the latter case, the curves can often be
assumed to form a simple random sample, in particular, the functional
observations are independent. For curves obtained from splitting a
continuous (in principle) time records into, say, daily or annual
curves, the assumption of independence is often violated.  This paper
focuses on the methodology and theory for the estimation of the mean
function of a functional time series, and on inference for the mean of
two functional time series. Despite their central importance, these
issues have not yet been studied.  The contribution of this
paper is thus two--fold: 1) we develop a methodology and
an asymptotic theory
for the estimation of the variance of the sample mean of temporally
dependent curves under model--free assumptions; 2) we propose
procedures for testing equality of two mean functions in
functional samples exhibiting temporal dependence.

\begin{figure} [!ht]
\begin{center}
\begin{capf} \label{f:rd}
{\small The horizontal component of the magnetic field measured in
one minute resolution at Honolulu magnetic observatory from 1/1/2001
00:00 UT to  1/7/2001 24:00 UT.}
\end{capf}
\includegraphics[width=.8\textwidth,angle=0]{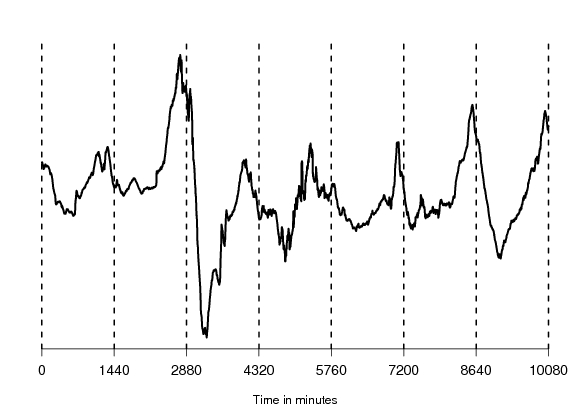}
\end{center}
\end{figure}

A functional time series $\{X_k,\,
k\in\mathbb{Z}\}$ is a sequence of curves $X_k(t)$, $t\in [a,b]$.
After normalizing to the unit interval, the curves are typically
defined as $X_k(t) = X(k+t), \ 0\le t \le 1$, where $\{X(u), u
\in\mathbb{R} \}$ is a continuous time record, which is often observed
at equispaced dense discrete points.  An example is given in
Figure~\ref{f:rd}, which shows seven consecutive functional
observations. More examples are studied in \citetext{HKbook}.
A central issue in the analysis of such data is to take into account
the temporal dependence of the observations.  The monograph of
\citetext{bosq:2000} studies the theory of linear functional time
series,  focusing on the functional
autoregressive model.  For many functional time series it is however
not clear what specific model they follow, and for many statistical
procedures it is not necessary to assume a specific model. In this
paper, we assume that the functional time series is stationary, but we
do not impose any specific model on it. We assume that the curves are
dependent in a very broad sense, which is made precise in
Section~\ref{ss:lpm}.
The dependence condition we use is however satisfied by all
models for functional time series used to date, including the linear,
multiplicative, bilinear and ARCH type processes.  We refer to
\citetext{hormann:kokoszka:2010} and
\citetext{aue:hormann:horvath:huskova:steinebach:2011} for examples.

A direct motivation for the research presented in this paper comes
from a two sample problem in which we wish to tests if the mean
functions of two functional time series are equal. A specific
problem, studied in greater depth in Section \ref{s:sd}, is to test if
the mean curves of certain financial assets are equal over certain periods.
This in turn allows us to conclude whether
the expectations of future market conditions are the same
or different at specific time periods.
In general, if the same mean is assumed for
the whole time series, whereas, in fact, it is different for disjoint
segments, the inference or exploratory analysis that follows will be
faulty, as all prediction and model fitting procedures for functional
time series start with subtracting the sample mean, viewed as an
estimate of the unique population mean function.
The same holds true for independent curves; if two subsamples have
different mean functions, subtracting the sample mean function
based on the whole data set will lead to spurious results.
Despite the
importance of two sample problems for functional data, they have
received little attention.  Recent papers of
\citetext{horvath:kokoszka:reimherr:2009} and
\citetext{panaretos:2010} are the only
contributions to a two sample problem in a functional setting which
develop inferential methodology.  \citetext{horvath:kokoszka:reimherr:2009}
compare linear operators in two functional regression models.
\citetext{panaretos:2010} focus
on testing the equality of the covariance operators in two samples of
iid Gaussian functional observations; our paper focuses on the means
of dependent (and non--Gaussian) observations.  We develop the
required methodology, justify it by asymptotic arguments, and describe
its practical implementation.

Any inference involving mean functions requires estimates of the
variability of the sample mean. In iid functional samples, the sample
covariance operator is used, but for functional time series this
problem is much more difficult. For scalar and vector--valued time
series, the variance of the sample mean is asymptotically approximated
by the long--run variance whose estimation has been one of the central
problems of time series analysis, studied in textbooks, see e.g.
\citetext{anderson:1971}, \citetext{brockwell:davis:1991},
\citetext{hamilton:1994}, and dozens of influential papers, see
\citetext{newey:west:1987}, \citetext{andrews:1991} and
\citetext{andrews:monahan:1992}, to name just a few.  Convergence of
various estimators of the long--run variance has been established under
several types of assumptions, including broad model specifications
(e.g. linear processes), cumulant conditions, and various mixing
conditions.  \citetext{hormann:kokoszka:2010} and
\citetext{gabrys:horvath:kokoszka:2010} advocated using the notion of
$L^p$--$m$--approximability for functional time series, as this
condition is intuitively appealing and is easy to verify for
functional time series models.  We therefore develop a general
framework for the estimation of the long--run covariance kernel in
this setting.  The long--run covariance kernel (equivalently operator)
is an infinite dimensional object, precise definitions are given in
Section~\ref{s:na}, whose estimation has not been studied yet.
\citetext{hormann:kokoszka:2010} studied the estimation of a long--run
covariance matrix obtained by projecting this operator onto a finite
dimensional basis. Their approach is based on the results for
vector--valued time series, and they use cumulant--like assumptions
which are difficult to verify for nonlinear functional time series.

The long--run covariance kernel corresponds to the asymptotic
variance in a normal approximation for the sample mean of a scalar
time series, but  no central limit theorem for a general
 functional time series has been established yet (results for
linear processes are established in \citetext{bosq:2000}).
We provide such a generally applicable result as well
(Theorem~\ref{t:2.1c}).

The remainder of the paper is organized as follows. We conclude the
introduction by defining the notion of dependence for functional time
series which we use throughout this paper.  Then, in
Section~\ref{s:na}, we state the asymptotic results for the mean of a
single functional time series, with proofs developed in Section
\ref{s:p-na}.  Section~\ref{s:te} focuses on the problem of testing
the equality of means of two functional samples which exhibit temporal
dependence.  In Section~\ref{s:sd}, we evaluate the finite sample
performance of the procedures proposed in Section~\ref{s:te} by means
of a simulation study and application to real data. The proofs of the
theorems stated in Section~\ref{s:te} are collected in
Section~\ref{s:p-te}.

\subsection{Approximable
functional time series} \label{ss:lpm}
We consider a stationary functional time series
$\{ \eg_i(t), \ i \in\mathbb{Z}, \ t\in [0,1] \}$, which we
can view as an error sequence in a more complex functional model,
for example a regression model, as in \citetext{gabrys:horvath:kokoszka:2010}.
We assume that these errors  are nonlinear moving averages
$\eg_i = f(\dg_i, \dg_{i-1}, \ldots)$,
for some measurable function $f: S^\infty \to L^2$, and iid elements $\dg_i$
of a measurable space $S$. In all models used in practice $S=L^2$.
To motivate the construction below, it is useful to write the $\eg_i$ as
\beq \label{e:1.6}
\eg_i = f(\dg_i, \ldots, \dg_{i-m+1}, \dg_{i-m},  \dg_{i-m-1} \ldots).
\eeq
Under (\ref{e:1.6}), the sequence $\{ \eg_i \}$
is stationary and ergodic. The function $f$
must decay sufficiently fast to ensure that the sequence $\{ \eg_i \}$ is
weakly dependent. The weak dependence condition is stated in terms
of an approximation by $m$--dependent sequences, namely, we require that
\beq \label{e:1.9}
\sum_{m\geq 1}
\lb E \int \lp \eg_i(t) - \eg_{i,m}(t) \rp^2 dt \rb^{1/2}<\infty,
\eeq
where
\beq \label{e:1.8a}
\eg_{i,m} = f(\dg_i, \ldots, \dg_{i-m+1}, \dg^{(m)}_{i, i-m},
\dg^{(m)}_{i, i-m-1}, \ldots ),
\eeq
with the sequences $\{ \dg_{i,k}^{(m)} \}$ being independent copies
of the sequence $\{ \dg_i \}$.
Note that the sum in (\ref{e:1.9}) does not depend on $i$.

The idea behind the above construction is that the function
$f$ decays so fast that the effect of the innovations $\dg_i$
far back in the past becomes negligible; they can be replaced by
different, fully independent innovations.
If the $\eg_i$ follow
a linear model $\eg_i = \sum_{j\ge 0} c_j(\dg_{i-j})$, condition
(\ref{e:1.9}) intuitively means that the approximations
by the finite moving averages $\eg_{i,m} = \sum_{0\le j\le m} c_j(\dg_{i-j})$
become increasingly precise. This means that the operators $c_j$ must decay
sufficiently fast in an appropriate operator norm.
We refer to \citetext{hormann:kokoszka:2010} and
\citetext{aue:hormann:horvath:huskova:steinebach:2011} for the details
and examples of nonlinear functional time series satisfying
(\ref{e:1.9}).

We also note that the general idea of using nonlinear moving averages
(Bernoulli shifts) and imposing moment conditions to quantify
dependence has been recently used in other contexts, see Wu
(\citeyear{wu:2005}, \citeyear{wu:2007}). The connections between such
notions and the traditional mixing conditions, or other notions of
weak dependence, e.g. that introduced by
\citetext{doukhan:louhichi:1999}, are only partially understood at
present. In particular, it is not clear which functional time series
models satisfy dependence conditions other than the approximability
(\ref{e:1.9}).

\section{Normal approximation
and long--run variance for
functional time series}
\label{s:na}

In this section, we state the central limit theorem for the sample
mean of an $L^2$--$m$--approximable functional time series.
Its applicability depends on the estimation of the covariance kernel
of the limit. We therefore also establish the consistency of the kernel
estimator of the long--run covariance kernel. We assume that
$\{ \eg_i \}$ is an $L^2$--$m$--approximable (and hence stationary)
functional time series satisfying
\beq \label{e:1.4}
E\eg_0 = 0, \ \  {\rm in}\ L^2
\eeq
and
\beq \label{e:1.8}
 \int E \eg_0^2(t) dt < \infty.
\eeq

\begin{Theorem}\label{t:2.1c} If
(\ref{e:1.6}), (\ref{e:1.9}), (\ref{e:1.4}), (\ref{e:1.8})  hold, then
\beq \label{co-1}
N^{-1/2}\sum_{i=1}^N\varepsilon_i\convd Z\;\;\;\mbox{in}\;\;\;L^2,
\eeq
where $Z$ is a Gaussian process with
\bd
EZ(t)=0 \ \ \ {\rm  and} \ \ \ E[Z(t)Z(s)] =c(t,s);
\ed
\beq \label{e:cts}
c(t,s) = E\eg_0(t)\eg_0(s) + \sum_{i \ge 1} E \eg_0(t)\eg_i(s)
+ \sum_{i \ge 1} E \eg_0(s)\eg_i(t).
\eeq
The infinite sums in the definition of the kernel $c$
converge in $L^2([0,1]\times[0,1])$, i.e.\ $c$ is a square
integrable function on the unit square.
\end{Theorem}

Theorem \ref{t:2.1c} is proven in Section \ref{s:p-na}.

The kernel $c$ is
defined analogously to the long--run variance of a scalar time series.
It is directly related to
the covariance operator of the sample mean defined by
\begin{align*}
\hat C_N(x) &= N
E \lb \lip \frac{1}{N} \sum_{i=1}^N\eg_i, x \rip
\frac{1}{N}\sum_{j=1}^N\eg_j  \rb\\
     &=\frac{1}{N}\sum_{i,j=1}^N E\lb \lip \eg_i, x\rip \eg_j  \rb.
\end{align*}
If the $\eg_i$ are independent, then
$\hat C_N(x) = N^{-1} \sum_{i=1}^N E \lb \lip \eg_i, x\rip \eg_i  \rb$
becomes the usual sample (empirical)
covariance operator, which plays a central role in many exploratory
and inferential tools of functional data analysis of iid functional
observations, mostly through the
empirical functional principal components defined as its eigenfunctions.
For functional time series, it is  not suitable. For any stationary
functional time series  $\{ \eg_i \}$,
\begin{align*}
\hat C_N(x)(t)
&= \int \lp\frac{1}{N}\sum_{i,j=1}^N E[\eg_i(s)\eg_j(t)] \rp x(s) ds\\
&= \int  c_N(t,s) x(s) ds,
\end{align*}
where
\beq \label{e:cN}
c_N(t,s) = \sum_{|k|<N} \lp 1- \frac{|k|}{N}\rp E [\eg_0(s)\eg_k(t)].
\eeq
The summands in (\ref{e:cN}) converge to those in (\ref{e:cts}),
but the estimation of the long--run covariance kernel $c$ is far
from trivial.

To enhance the
applicability of our result, we state it for the case of a nonzero mean
function, which is estimated by the sample mean. We thus assume that
\beq \label{e:sample1def}
X_i(t) = \mu(t) + \eg_i(t), \ \ \ 1 \le i \le N,
\eeq
with the series $\{ \eg_i \}$ satisfying the assumptions of
Theorem~\ref{t:2.1c}.

Let $K$ be a kernel (weight) function defined on the line and satisfying
the following conditions:
\beq \label{e:3.1}
K(0) =1,
\eeq
\beq \label{e:3.2}
K \ {\rm is \ continuous},
\eeq
\beq \label{e:3.3}
K \ {\rm is \ bounded},
\eeq
\beq \label{e:3.4}
K(u) = 0,\ {\rm if }\ |u| > c, \ \     {\rm for\ some}\ c>0.
\eeq
Condition (\ref{e:3.4}) is assumed only to simplify the proofs,
a sufficiently fast decay could be assumed instead.

Next we define the empirical (sample) correlation functions
\beq \label{e:gamma}
\hat\ga_i(t,s) = \frac{1}{N} \sum_{j=i+1}^N
\lp X_j(t) - \bar X_N(t) \rp\lp X_{j-i}(s) - \bar X_N(s) \rp,
\ \ \ 0 \le i \le N-1,
\eeq
where
\bd
\bar X_N(t) = \frac{1}{N} \sum_{1 \le i \le N} X_i(t).
\ed
The estimator for $c$ is  given by
\beq \label{e:hatC}
\hat c_N(t,s) = \hat\ga_0(t,s)
+ \sum_{i=1}^{N-1} K\lp \frac{i}{h} \rp
\lp \hat\ga_i(t,s)+\hat\ga_i(s,t)  \rp
\eeq
where $h= h(N)$ is the  smoothing bandwidth satisfying
\beq \label{e:3.5}
h(N) \to \infty, \ \ \frac{h(N)}{N} \to 0, \ \ \ {\rm as} \ N\to \infty.
\eeq
In addition to (\ref{e:1.9}), we also assume that
\beq \label{m1}
\lim_{m\rightarrow \infty}
m\lb E \int \lp \eg_n(t) - \eg_{n,m}(t) \rp^2 dt \rb^{1/2}=0.
\eeq

\begin{Theorem} \label{t:3.1c} Suppose the functional
time series $\{ X_i \}$ follows model (\ref{e:sample1def}).
Under conditions
(\ref{e:1.6}), (\ref{e:1.9}), (\ref{e:1.4}), (\ref{e:1.8}),
(\ref{e:3.1})--(\ref{e:3.4}), (\ref{e:3.5}), (\ref{m1}),
\beq \label{e:6.1}
\dint \lp \hat c_N(t,s) - c(t,s) \rp^2dtds\convP 0,
\eeq
with $c(t,s)$ by defined by (\ref{e:cts}) and
$\hat c_N(t,s)$ by (\ref{e:hatC}).
\end{Theorem}

Theorem \ref{t:3.1c} is proven in Section \ref{s:p-na}. First
we use  the results of this section in the problem of testing the
equality of means in two functional samples.

\section{Testing the equality of
mean functions} \label{s:te}

We consider two samples of curves, $X_1, X_2, \ldots, X_N$ and
$X_1^*, X_2^*, \ldots X_M^*$,
satisfying the following location models
\beq \label{e:sample2def}
X_i(t) = \mu(t) + \eg_i(t),
\ \ \
X_j^*(t) = \mu^*(t) + \eg_j^*(t).
\eeq
The error functions $\eg_i$ are assumed to satisfy the conditions
stated in Sections \ref{s:I} and \ref{s:na}. The functions $\eg_j^*$
are assumed to satisfy exactly the same conditions. In particular,
their long--run covariance kernel is defined by
\bd
c^*(t,s) = E\eg_0^*(t)\eg_0^*(s) + \sum_{i \ge 1} E \eg_0^*(t)\eg_i^*(s)
+ \sum_{i \ge 1} E \eg_0^*(s)\eg_i^*(t).
\ed
We assume that
\beq \label{e:1.1}
\lbr \eg_i, \ 1 \le i \le N \rbr \ \ {\rm and} \ \
\lbr \eg_j^*, \ 1 \le j \le M \rbr \ \ \ {\rm are \ independent}.
\eeq

We are interested in testing
\beq \label{e:1.2}
H_0:\ \mu = \mu^*
\eeq
against the alternative
\beq \label{e:1.3}
H_A:\ \mu \neq \mu^*.
\eeq
The equality in (\ref{e:1.2}) is in the space $L^2=L^2([0,1])$, i.e.
$\mu=\mu^*$ means that
$\int \lp \mu(t) - \mu^*(t) \rp^2 dt = 0$, and  the
alternative means that  $\int \lp \mu(t) - \mu^*(t) \rp^2 dt > 0$.

Since the statistical inference is about the mean functions of the
observations, our procedures are based on the sample mean curves
\bd
\bar X_N(t) = \frac{1}{N} \sum_{1 \le i \le N} X_i(t)
\ \ \ {\rm and}\ \ \
\bar X_M^* =  \frac{1}{M} \sum_{1 \le i \le M} X_j(t).
\ed
The sample means  $\bar X_N$ and $\bar X_M^*$ are unbiased estimators
of $\mu$ and $\mu^*$, respectively, so $H_0$ will be rejected if
\bd
U_{N,M} = \frac{NM}{N+M} \int \lp\bar X_N(t)-\bar X_M^*(t)  \rp^2 dt
\ed
is large.


Before introducing the test procedures, we state two results which
describe the asymptotic behavior of the statistic $U_{N,M}$ under
$H_0$ and $H_A$. They  motivate and explain the development that
follows.

\begin{Theorem} \label{t:1.1}
Suppose $H_0$, the assumptions of Theorem \ref{t:2.1c} (and analogous
assumptions for the $\eg_j^*$) and (\ref{e:1.1}) hold.
If
\beq \label{e:1.9a}
\frac{N}{N+M} \to \theta, \ \ \ {\rm for \ some} \ \ 0\le  \theta \le 1, \ \ \
{\rm as} \ \min(M,N) \to \infty,
\eeq
then
\bd
U_{N,M} \convd \int_0^1 \Gg^2(t) dt,
\ed
where $\lbr \Gg(t), \ 0 \le t \le 1 \rbr$ is a mean zero  Gaussian process
with covariances
\bd
E[\Gg(t)\Gg(s)] =
d(t,s) := (1-\theta) c(t,s) + \theta c^*(t,s).
\ed
\end{Theorem}

\begin{Theorem} \label{t:1.2}
If  $H_A$, and the remaining assumptions of Theorem~\ref{t:1.1} hold,
then
\bd
\frac{N+M}{NM} U_{N, M} \convP \int_0^1 \lp \mu(t) - \mu^*(t) \rp^2 dt.
\ed
In particular, if $0< \theta < 1$, then $ U_{N, M} \convP \infty$.
\end{Theorem}

The kernel $d(t,s)$ in Theorem \ref{t:1.1} defines a covariance operator
$D$. The eigenvalues of $D$ are nonnegative, and are denoted by
$\la_1 \ge \la_2 \ge \ldots$. By the Karhunen--Lo{\'e}ve expansion, we
have
\beq \label{e:1.10}
\int_0^1 \Gg^2(t) dt = \sum_{i=1}^\infty \la_i N_i^2,
\eeq
where $\{ N_i, \ 1 \le i < \infty \}$ are independent standard normal
random variables.

Since the eigenvalues $\la_i$ are unknown, the right--hand side of
(\ref{e:1.10}) cannot be used directly to simulate the distribution
of $\int \Gg^2(t) dt $. We will now explain how to estimate the $\la_i$'s.

Suppose $\hat D_{N,M}$ is an $L^2$--consistent estimator of $D$ , i.e.
\beq \label{e:1.11}
\dint \lp\hat d_{N,M}(t,s) - d(t,s) \rp^2 dtds \convP 0, \ \ \
{\rm as} \ \min(M,N) \to \infty.
\eeq

We discuss the construction of estimators $\hat D_{N,M}$ satisfying
(\ref{e:1.11}) below. For the estimators we propose, relation
(\ref{e:1.11}) holds regardless whether $H_0$ or $H_A$ holds,
they do not depend on $\mu$ or $\mu^*$ either. We will also see
that the critical relations (\ref{e:1.13})
hold under $H_A$ as well as under $H_0$. The distribution of
$\int \Gg^2(t) dt$ can thus be estimated also under the alternative.

Let
\bd
\hat\la_1 = \hat\la_1(N,M) \ge \hat\la_2 = \hat\la_2(N,M)\ge \ldots
\ed
denote the eigenvalues of $\hat D_{N,M}$, i.e.
\beq \label{e:11a}
\int \hat d_{N,M}(t,s) \hat\fg_i(s)ds = \hat\la_i  \hat\fg_i(t),
\eeq
where the $\hat\fg_i(t) = \hat\fg_i(t; N, M)$ are the corresponding
eigenfunctions satisfying $\int \hat\fg_i^2(t)dt =1$.
Choosing $p$ so large that $\sum_{i=1}^p \hat\la_i$ is
a large percentage of  $\sum_{i=1}^{N+M} \hat\la_i$, we can approximate
the distribution of $\int \Gg^2(t) dt$ by that of
$\sum_{i=1}^p \hat\la_i N_i^2$.

The statistical inference is based on the difference
$\bar X_N - \bar X_M^*$. Observe that
\bd
\frac{MN}{M+N}
E \lb \lp \bar X_N(t) - \bar X_M^*(t)\rp
 \lp \bar X_N(s) - \bar X_M^*(s)\rp \rb  \to d(t,s), \ \ \
{\rm as} \ \min(M,N) \to \infty,
\ed
that is, $d$ is the asymptotic covariance kernel of the difference
$\bar X_N - \bar X_M^*$.
We therefore  use projections onto the eigenfunctions
$\fg_1, \fg_2, \ldots, \fg_p$ associated with the $p$ largest
eigenvalues of $D$. This is analogous to projecting onto
 the functional principal components in one sample problems,
as these form an $L^2$--optimal orthonormal basis.
Without any loss of generality, we
assume that the $\fg_1, \fg_2, \ldots, \fg_p$
form an orthonormal system (the $\fg_i$ are orthogonal under
(\ref{e:1.12}), so only a normalization to unit norm is required).
We define the projections
\beq \label{e:2.1}
a_i = \lip \bar X_N - \bar X_M^*, \fg_i \rip, \ \ 1 \le i \le p,
\eeq
and the vectors
\bd
\ba = [a_1, a_2, \ldots, a_p]^T.
\ed
We show in the proof of Theorem \ref{t:2.1} that
\beq \label{e:2.2}
\lp \frac{MN}{M+N}\rp^{1/2} \ba \convd {\bf N}_p(\bzero, \bQ),
\eeq
where ${\bf N}_p(\bzero, \bQ)$ stands for the $p$--variate normal
random vector with mean zero and the covariance matrix
 $\bQ={\rm diag}(\la_1, \la_2, \ldots, \la_p)$.
Since the operator $D$ is unknown, we cannot compute the $\fg_i$.
However, any estimator for $D$ satisfying (\ref{e:1.11}) can be used
to find estimates for the $\fg_i$. Let $\hat\fg_i$ be the empirical
eigenfunctions defined by (\ref{e:11a}), and set
\bd
\hat a_i = \lip \bar X_N - \bar X_M^*, \hat\fg_i \rip, \ \ 1 \le i \le p.
\ed
The limit relation (\ref{e:2.2}) suggests the following statistics:
\beq \label{e:stat1}
U_{N,M}^{(1)} =  \frac{MN}{M+N} \sum_{i=1}^p {\hat a_i}^2
\eeq
and
\beq \label{e:stat2}
U_{N,M}^{(2)} =  \frac{MN}{M+N} \sum_{i=1}^p \frac{{\hat a_i}^2}{\hat\la_i}.
\eeq

The following theorem
establishes the limits of $U_{N,M}^{(1)}$ and $U_{N,M}^{(1)}$ under $H_0$.

\begin{Theorem} \label{t:2.1}
Suppose $H_0$, the remaining assumptions of Theorem \ref{t:1.1},
(\ref{e:1.11}) and
\beq \label{e:1.12}
\la_1 > \la_2 > \ldots > \la_p > \la_{p+1}
\eeq
hold. Then
\beq \label{e:2.3}
U_{N,M}^{(1)} \convd \sum_{i=1}^p \la_i N_i^2,
\eeq
where $N_1, N_2, \ldots, N_p$ are independent standard normal
random variables, and
\beq \label{e:2.4}
U_{N,M}^{(2)} \convd \chi^2(p),
\eeq
where $\chi^2(p)$ stands for a chi-square random variable with $p$ degrees
of freedom.
\end{Theorem}

We note that $U_{N,M}^{(1)}$ is essentially the first $p$ terms in
the Karhunen--Lo{\'e}ve expansion of the integral in the definition of
 $U_{N,M}$. Thus, the limit in (\ref{e:2.3}) is exactly the random
variable we used to approximate the distribution of  $U_{N,M}$.
The limit in (\ref{e:2.4}) is distribution free.

The consistency of the tests based on $U_{N,M}^{(1)}$ and $U_{N,M}^{(2)}$
follows from the following result.

\begin{Theorem} \label{t:2.2}
Suppose $H_A$, the remaining assumptions of Theorem \ref{t:1.1},
(\ref{e:1.11}) and (\ref{e:1.12}) hold. Then
\bd
\frac{N+M}{NM} U_{N,M}^{(1)}
\convP \sum_{i=1}^p \lip \mu - \mu^*, \fg_i \rip^2
\ed
and
\bd
\frac{N+M}{NM} U_{N,M}^{(2)}
\convP \sum_{i=1}^p \frac{\lip \mu - \mu^*, \fg_i \rip^2}{\la_i}.
\ed
In particular, if $0<\theta < 1$ in (\ref{e:1.9a}), then
$U_{N,M}^{(1)}\convP \infty$ and $U_{N,M}^{(2)}\convP \infty$,
provided $\lip \mu - \mu^*, \fg_i \rip \neq 0$
for at least one $1 \le i \le p$.
\end{Theorem}

We see that the condition for the consistency is that $\mu - \mu^*$
is not orthogonal to the linear subspace of $L^2$ spanned by the
eigenfunctions $\fg_i, \ 1 \le i \le p.$

\vspace{2mm}

The implementation of the tests based on Theorems \ref{t:2.1} and \ref{t:2.2}
depends on the existence of an estimator of the kernel $d(t,s)$ which
satisfies (\ref{e:1.11}). The remainder of this section is dedicated
to this issue.

The estimation of $D$ is very simple if the $\eg_i$ are iid, and the
$\eg_j^*$ are iid. In this case, setting,
\bd
\hat\theta = \frac{N}{N+M},
\ed
we can use
\beq \label{e:hatD}
\hat d_{N,M}(t,s) = (1-\hat\theta) \hat c_N(t,s)
+ \hat\theta \hat c_M^*(t,s),
\eeq
where
\bd
\hat c_N(t,s) = \frac{1}{N} \sum_{i=1}^N
 \lp X_i(t) - \bar X_N(t)  \rp \lp X_i(s) - \bar X_N(s)  \rp;
\ed
\bd
\hat c_M^*(t,s) = \frac{1}{M} \sum_{j=1}^M
 \lp X_j^*(t) - \bar X_M^*(t)  \rp \lp X_j^*(s) - \bar X_N^*(s)  \rp.
\ed
By condition (\ref{e:1.8}), we can use the weak law of large numbers
in a Hilbert space to establish (\ref{e:1.11}).
The estimation of $D$ is much more difficult if only (\ref{e:1.9}) is
assumed, and its asymptotic justification relies on Theorem~\ref{t:3.1c}.
Recall the definition of the estimator $\hat c_N(t,s)$ given in (\ref{e:hatC}),
and define the estimator $\hat c_M^*(t,s)$ fully analogously.
Our estimator for $d(t,s)$ is then (\ref{e:hatD}) with
$\hat c_N(t,s)$ and $\hat c_M^*(t,s)$ so defined. The following result
then follows directly from Theorem~\ref{t:3.1c}.

\begin{Theorem} \label{t:3.1}
Suppose  the functional time series $\{ X_i \}$ satisfies the assumptions
of Theorem~\ref{t:3.1c}, and the series $\{ X_j^* \}$ satisfies the
same  assumptions stated in terms of the $\eg_j^*$. If (\ref{e:1.1}) holds,
then  (\ref{e:1.11}) holds.
\end{Theorem}

We emphasize that under the conditions of Theorem \ref{t:3.1}
relation (\ref{e:1.11}) holds both under $H_0$ and  $H_A$.

\vspace{2mm}

We now focus on the numerical issues related to the computation of the
$\hat a_i$ and the $\hla_i$ appearing in the definitions of statistics
$U_{N, M}^{(1)}$ and $U_{N, M}^{(2)}$. The $\hat a_i$ and the $\hla_i$
require the computation of the eigenfunctions and the eigenvalues of
the operator $\hat D$. Except in the case of independent observations
in each of the two samples, these quantities cannot be computed using
existing software because   $\hat D$ is not an empirical covariance
operator of a functional iid sample. We recommend the following
algorithm which we used to implement the tests.
Let $\{ e_{\ell}, \ell \ge 1 \}$ be an orthonormal basis.
The results reported in Section~\ref{s:sd} are based on an implementation
which uses the Fourier basis. In order to find approximate solutions to
    \begin{equation}\label{e:eigen}
        \int \hat{d}_{N,M}(t,s) \phi(s) ds = \lambda \phi(t),
    \end{equation}
    we use approximate expansions for $\phi(s)$ and $\hat{d}_{N,M}(t,s)$:
    $$
    \begin{aligned}
    \phi(s) &\approx \sum_{k=1}^{49} \phi_k e_k(t),\\
    \hat{d}_{N,M}(t,s) &\approx \sum_{k=1}^{49} \sum_{\ell=1}^{49} d_{k,\ell} e_k(t) e_{\ell}(s).\\
    \end{aligned}
    $$

    The coefficients $\phi_k$ and $d_{k,\ell}$ are given by
    $$
        \phi_k = \int \phi(t) e_k(t) dt
    $$
    and
    \beq \label{e:coefD}
        d_{k,\ell} =\int \int \hat{d}_{N,M}(t,s) e_k(t) e_{\ell}(s) dtds.
    \eeq
    We replace $\phi(s)$ and $\hat{d}_{N,M}(t,s)$ in the left side of \eqref{e:eigen} with the above expansions to obtain
    \begin{equation}\label{e:step1}
        \sum_{k=1}^{49} \sum_{\ell=1}^{49} d_{k,\ell} e_k(t) \approx \lambda \phi(t).
    \end{equation}
    Multiplying both sides of \eqref{e:step1} by $e_j$ and integrating yields
    $$
        \sum_{\ell=1}^{49} d_{j,\ell} \phi_{\ell} \approx \lambda \phi_j,\ \ \ 1\le j \le 49.
    $$
    In matrix form this is
    $$
    {\bf D}{\boldsymbol \phi} = \lambda {\boldsymbol \phi},
    $$
    where ${\bf D}=[d_{j,\ell},\ 1\le j, \ell \le 49 ]$ and ${\boldsymbol \phi}=[\phi_1, \phi_2, \dots, \phi_{49}]^T$. Thus we have reduced the problem of finding solutions of \eqref{e:eigen} to finding eigenvalues and eigenvectors of the matrix ${\bf D}$.  Let $[\phi_{m,1}, \phi_{m,2}, \dots, \phi_{m,49}]^T$
    be the eigenvector corresponding to the $m$th
    largest eigenvalue of ${\bf D}$.
    Then the eigenfunction associated with the $m$th largest eigenvalue of $\hat{d}_{N,M}(t,s)$ is approximately $\sum_{\ell}^{49} \phi_{m,\ell} e_{\ell}(t)$.  Using this notation, we obtain
    \beq \label{e:hata}
    \begin{aligned}
        \hat{a}_i &= \lip \bar X_N - \bar X_M^*, \hat\fg_i \rip\\
        &\approx \sum_{\ell=1}^{49} \phi_{m,\ell} \left(N^{-1}\sum_{i=1}^N \lip X_i, e_{\ell} \rip - M^{-1} \sum_{j=1}^M \lip X_j^*, e_{\ell} \rip \right).\\
    \end{aligned}
    \eeq

We then obtain $U_{N,M}^{(1)}$ and $U_{N,M}^{(2)}$
as per \eqref{e:stat1} and \eqref{e:stat2}.

\vspace{2mm}

To complete the description of the test procedures, we must specify how
the value of $p$ in the definitions of  $U_{N,M}^{(1)}$ and $U_{N,M}^{(2)}$
is selected. This issue has been extensively studied in one sample
problems, and several approaches have been put forward, including
cross validation and penalty criteria. In our experience (for smooth
densely recorded curves), the simple cumulative
variance method advocated by \citetext{ramsay:silverman:2005} has been
satisfactory. We therefore recommend to use $p$ such that the first $p$
empirical functional principal components in each sample explain about
85\% of the variance. As we will see in the data examples studied in
Section~\ref{s:sd}, it is typically useful to look at the P--values
for a range of $p$'s.

\section{A simulation study
and data examples}\label{s:sd}

We begin by  presenting the results of a simulation study intended to evaluate
the empirical size and power of the testing procedures introduced in
Section \ref{s:te}. We then illustrate their properties on two data
examples.

\subsection{A simulation study} \label{ss:sim}
In this section, we compare the the performance of the tests based on
statistics $U_{N,M}^{(1)}$ and $U_{N,M}^{(2)}$ using simulated Gaussian
functional data. We consider all combinations of sample sizes
$N,M=50, 100, 100$, and each pair of data generated processes
was replicated three thousand times.
To investigate the empirical
size, without loss of generality, we set $\mu(t) = \mu^*(t) = 0$.
Under the alternative, we set $\mu(t)=0$ and $\mu^*(t) = at(1-t)$.
The power is then a function of the parameter $a$.  We considered two
settings for the errors:
\begin{enumerate}
\item   Both the $\eg_i(t)$ and the $\eg_j^*(t)$ are iid
Brownian bridges.
\item   Both the $\eg_i(t)$ and the $\eg_j^*(t)$ are
functional AR(1) (FAR(1)) processes with the kernel,
$$
\psi(t,s)=\frac{e^{-(t^2+s^2)/2}}{4\int_0^1 e^{-x^2} dx}.
$$
That is, the error terms, $\eg_i(t)$, follow the model
$$
\eg_i(t) = \int_0^1 \psi(t,s) \eg_{i-1}(s) ds + B_i(t),
$$
where $B_i(t)$ are iid Brownian bridges.
\end{enumerate}

We calculated the test statistics  $U_{N,M}^{(1)}$ and $U_{N,M}^{(2)}$
as explained in Section~\ref{s:te}.
These statistics depend on the choice of the weight functions $K$ and
$K^*$, and the bandwidth functions $h$ and $h^*$. A great deal of attention
has been devoted over several decades
to the optimal selection of these functions for
scalar and vector--valued time series, and we cannot address this issue
within the space of this paper.
We follow the recommendation of \citetext{politis:romano:1996}
and use, for both samples, the flat top kernel
\bd
        K(t)=   \begin{cases}
                    1     	& 0 \le |t| < 0.1,\\
                    1.1-|t|    & 0.1 \le |t| < 1.1,\\
                    0        & 1.1 \le |t|\\
                \end{cases}
\ed
with $h=N^{1/3}$ and  $h^*=M^{1/3}$.
We emphasize that this full estimation
procedure was used for all data generating processes, including those
with independent errors.

\begin{table}[ht!]
\begin{capt}\label{t:bridge}{\small Power of Test (in \%)
using $U_{100,200}^{(1)}$ and $U_{100,200}^{(2)}$
with iid Brownian-bridge errors}
\end{capt}
\centering
\begin{tabular}{|c|r|r|r|r|r|r|}

\hline

                                &             \multicolumn{2}{c|}{$\alpha=.01$}                                          &             \multicolumn{2}{c|}{$\alpha=.05$}                                                &            \multicolumn{2}{c|}{$\alpha=.10$}                                         \\ \hline

$              a              $&$        U_{100, 200}^{(1)}           $&$        U_{100, 200}^{(2)}           $&$        U_{100, 200}^{(1)}           $&$                U_{100, 200}^{(2)}           $&$        U_{100, 200}^{(1)}           $&$        U_{100, 200}^{(2)}          $\\ \hline

$              0.0          $&$        1.5          $&$        1.5          $&$        6.3          $&$        6.2          $&$        11.4        $&$        11.6        $\\ \hline

$              0.1          $&$        2.5          $&$        3              $&$        7.4          $&$        8.0          $&$        13.2        $&$        13.6        $\\ \hline

$              0.2          $&$        6.0          $&$        4.4          $&$        16.7        $&$        13.0        $&$        24.8        $&$        20.2        $\\ \hline

$              0.3          $&$        14.2        $&$        9.2          $&$        30.4        $&$        23.2        $&$        41.3        $&$        33.0        $\\ \hline

$              0.4          $&$        26.4        $&$        17.0        $&$        48.5        $&$        36.1        $&$        60.8        $&$        48.0        $\\ \hline

$              0.5          $&$        44.0        $&$        31.4        $&$        64.7        $&$        53.5        $&$        75.2        $&$        64.3        $\\ \hline

$              0.6          $&$        59.4        $&$        45.9        $&$        80.3        $&$        68.0        $&$        87.9        $&$        78.4        $\\ \hline

$              0.7          $&$        78.0        $&$        64.7        $&$        91.8        $&$        82.4        $&$        96.0$&$                89.0        $\\ \hline

$              0.8          $&$        88.0        $&$        78.0        $&$        95.9        $&$        90.9        $&$        98.0        $&$        94.6        $\\ \hline

$              0.9          $&$        94.2        $&$        88.1        $&$        98.7        $&$        95.8        $&$        99.4        $&$        97.9        $\\ \hline

$              1.0          $&$        98.0        $&$        94.8        $&$        99.5        $&$        98.5        $&$        99.9        $&$        99.3        $\\ \hline

$              1.1          $&$        99.6        $&$        98.4        $&$        100.0     $&$        99.8        $&$        100.0     $&$        100.0     $\\ \hline

$              1.2          $&$        99.9        $&$        99.4        $&$        100.0     $&$        99.9        $&$        100.0     $&$        100.0     $\\ \hline

$              1.3          $&$        100.0     $&$        99.8        $&$        100.0     $&$        100.0     $&$        100.0     $&$        100.0     $\\ \hline
\end{tabular}
\end{table}

\begin{table}[ht!]
\begin{capt}\label{t:far}{\small Power of Test (in \%)
using $U_{100,200}^{(1)}$ and $U_{100,200}^{(2)}$ with FAR(1) errors}
\end{capt}
\centering
\begin{tabular}{|c|r|r|r|r|r|r|}

\hline

                                &             \multicolumn{2}{c|}{$\alpha=.01$}                                          &             \multicolumn{2}{c|}{$\alpha=.05$}                                                &            \multicolumn{2}{c|}{$\alpha=.10$}                                         \\ \hline

$              a              $&$        U_{100, 200}^{(1)}           $&$        U_{100, 200}^{(2)}           $&$        U_{100, 200}^{(1)}           $&$                U_{100, 200}^{(2)}           $&$        U_{100, 200}^{(1)}           $&$        U_{100, 200}^{(2)}          $\\ \hline

$              0.0          $&$        1.8          $&$        1.9          $&$        6.6          $&$        7.2          $&$        12.2        $&$        13.5        $\\ \hline

$              0.1          $&$        2.4          $&$        2.2          $&$        7.9          $&$        7.7          $&$        13.5        $&$        14.5        $\\ \hline

$              0.2          $&$        5.1          $&$        3.3          $&$        13.6        $&$        11.6        $&$        21.6        $&$        18.7        $\\ \hline

$              0.3          $&$        9.8          $&$        6.3          $&$        23.6        $&$        17.6        $&$        34.5        $&$        26.8        $\\ \hline

$              0.4          $&$        19.4        $&$        12.3        $&$        35.9        $&$        26.5        $&$        46.7        $&$        36.3        $\\ \hline

$              0.5          $&$        26.8        $&$        19.5        $&$        47.9        $&$        38.6        $&$        60.4        $&$        49.7        $\\ \hline

$              0.6          $&$        42.1        $&$        29.6        $&$        62.2        $&$        51.8        $&$        73.1        $&$        62.5        $\\ \hline

$              0.7          $&$        56.4        $&$        42.8        $&$        75.4        $&$        63.8        $&$        83.2        $&$        74.0        $\\ \hline

$              0.8          $&$        68.6        $&$        53.8        $&$        85.7        $&$        74.6        $&$        91.5        $&$        83.1        $\\ \hline

$              0.9          $&$        80.8        $&$        67.6        $&$        92.7        $&$        85.9        $&$        96.4        $&$        91.9        $\\ \hline

$              1.0          $&$        87.4        $&$        78.7        $&$        95.9        $&$        90.8        $&$        98.1        $&$        94.5        $\\ \hline

$              1.1          $&$        93.7        $&$        86.8        $&$        97.9        $&$        95.8        $&$        99.1        $&$        97.6        $\\ \hline

$              1.2          $&$        97.6        $&$        93.7        $&$        99.5        $&$        98.1        $&$        99.8        $&$        99.2        $\\ \hline

$              1.3          $&$        98.5        $&$        96.4        $&$        99.7        $&$        98.9        $&$        99.9        $&$        99.6        $\\ \hline

\end{tabular}

\end{table}

\begin{figure}
\begin{center}
\begin{capf} \label{f:a0.8}
{\small Fifty trajectories of the Brownian bridge (left)
and fifty independent trajectories of the Brownian bridge
plus $\mu^*(t) = 0.8 t(1-t)$ (right). The tests can detect
the different means with probability close to 90\%.}
\end{capf}
\includegraphics[height=6cm,width=.45\textwidth,angle=0]{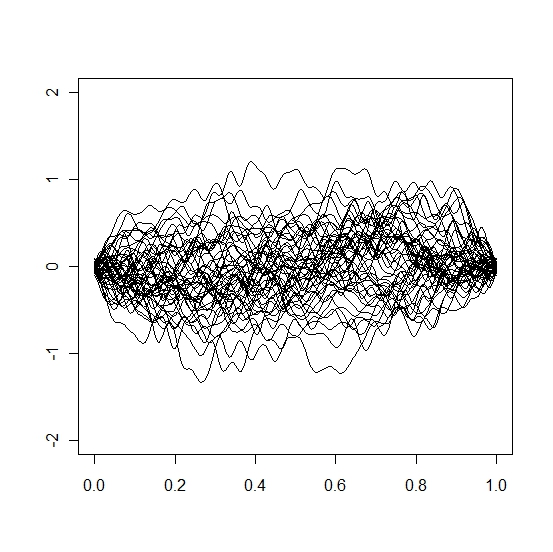}
\includegraphics[height=6cm,width=.45\textwidth,angle=0]{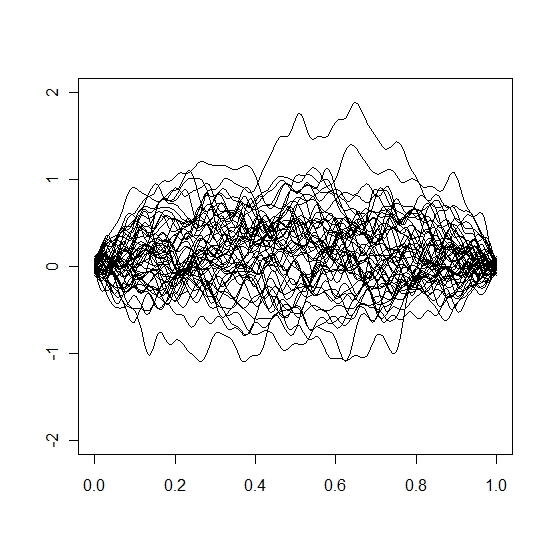}
\end{center}
\end{figure}

The results of the simulation study can be summarized as follows.  The
empirical size of the tests is larger in the case of FAR(1) errors.
When $a$ increases to $0.2$ or larger, the empirical power of the test
is smaller in the case of FAR(1) errors.  Thus increasing the
dependence in the error terms increases the size and decreases the
power of the test.  In both cases the tests have have a slightly
larger-than-nominal size and very good power.  These observations are
illustrated in Tables \ref{t:bridge} and \ref{t:far}. Based on the
whole simulation study, we can conclude that the performance of both
tests is better if the sample sizes $N$ and $M$ are about equal.
For example, for $N=M=100$, the empirical sizes are closer to the nominal
sizes then in the case $N=100, M=200$ shown in Tables
\ref{t:bridge} and \ref{t:far}. The
power is very high even for small sample sizes. This is illustrated in
Figure~\ref{f:a0.8} which shows the samples with $N=M=50$ and with
slightly different means ($a=0.8$). Visual inspection does not readily
lead to the conclusion that the samples in the left and right panels
of Figure~\ref{f:a0.8} have different means, yet our tests can detect
it with a very high probability.  None of the two test statistics
clearly dominates the other for the simulated Gaussian data, but a
difference in behavior can be see when the tests are applied to real
data sets, to which we now turn.

\subsection{Mediterranean fruit flies}\label{ss:medfly}
\begin{figure}
\begin{capf}\label{f:egg-curves}
{\small  Ten randomly selected smoothed egg-laying curves of short-lived
medflies (left panel), and ten such curves for long--lived medfies
(right panel).}
\end{capf}
\begin{center}
\includegraphics[height=8cm, width=7cm]{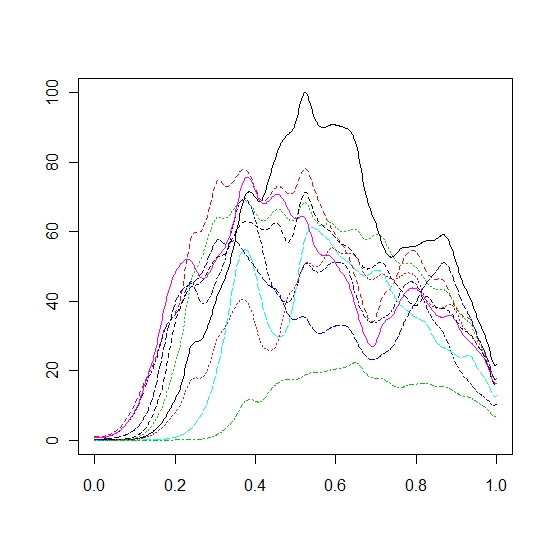}
\includegraphics[height=8cm, width=7cm]{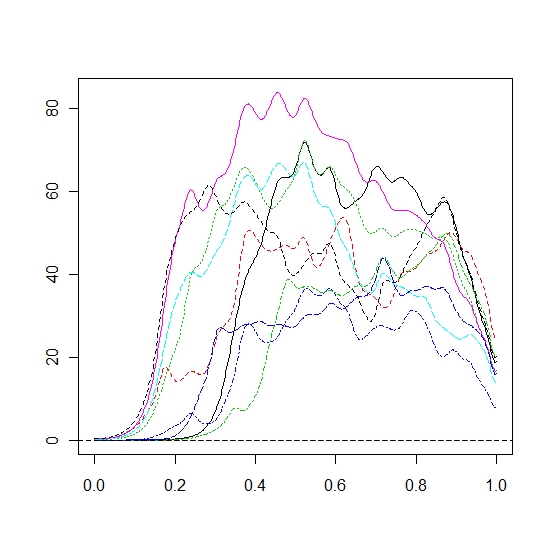}
\end{center}
\end{figure}

In our first example, it can be assumed that the curves in each
sample are independent, as they were obtained from a randomized experiment.
The data set used in this example was kindly made available to
us by Hans--Georg M{\"u}ller. It was extensively studied in biological
and statistical literature, see \citetext{muller:stadtmuller:2005}
and references therein.
We consider 534 egg-laying
curves (count of eggs per unit time interval)
of medflies who lived at least 30 days. Each function is
defined over an interval $[0,30]$, and its value on day $t\le 30$ is
the count of eggs laid by fly $i$ on that day.  The 534 flies are
classified into long-lived, i.e. those who lived longer than 44 days,
and short-lived, i.e. those who died before the end of the 44th day
after birth. In the data set, there are 256 short-lived, and 278
long-lived flies.
This classification naturally defines two samples:
{\em Sample 1:}
 the egg-laying curves
$\{X_i(t),\  0< t \le 30, \ i=1,2, \ldots, 256\}$ of the  short-lived
flies.
{\em Sample 2:}
 the egg-laying curves
$\{X_j^*(t),\  0< t \le 30, \ j=1,2, \ldots, 278\}$ of the  long-lived
flies.
The egg-laying curves are very irregular; Figure \ref{f:egg-curves}
shows ten smoothed curves of short- and long-lived flies.
The tests are applied to such smooth trajectories.

\begin{table}
\begin{capt} \label{tb:mean-flies}
{\small P-values (in percent) of the tests based on statistics
 $U_{N,M}^{(2)}$ and $U_{N,M}^{(1)}$ applied to medfly data.}
\end{capt}
\begin{center}
\begin{tabular}{|c|rrr rrr rrr|}
\hline
$p$      & 1  & 2 & 3 & 4 & 5 & 6 & 7 & 8 & 9 \\
\hline
$U^{(1)}$&1.0 &1.0&1.0&1.1 &1.1 &1.0 &1.0&1.1&1.1\\
$U^{(2)}$&1.0 &2.2&3.0&5.7 &10.3&15.3&3.2&2.7&5.0\\
\hline
\end{tabular}
\end{center}
\end{table}

\begin{figure}
\begin{center}
\begin{capf} \label{f:mean-flies}
{\small Estimated mean functions for the medfly data:\
short lived --solid line;\ long lived --dashed line.  }
\end{capf}
\includegraphics[height=7cm,width=.6\textwidth,angle=0]{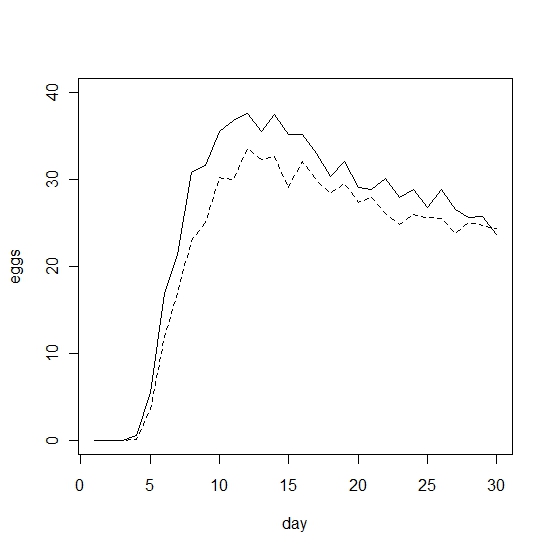}
\end{center}
\end{figure}

Table~\ref{tb:mean-flies} shows the P--values as a function of
$p$. For both samples, $p=2$ explains slightly over 85\%
of the variance, so this is the value we would recommend using.
Both tests reject the equality of the mean functions, even though
the sample means, shown in Figure \ref{f:mean-flies}, are not
far apart.
The P-values for the statistic $U^{(1)}$ are much more stable,
equal to about 1\%, no matter the value of $p$. The behavior
of the test based on $U^{(2)}$ is more erratic. This indicates
that while the test based on $U^{(2)}$ is easier to apply
because it uses standard chi-square critical values, the test
based on $U^{(1)}$ may be more reliable.

\subsection{Eurodollar futures contracts} \label{ss:eurodollar} Our
next example uses financial data kindly made available by Vladislav
Kargin.  This data is used as an example of modeling with the
functional AR(1) process in \citetext{kargin:onatski:2008}. The
curves, one curve per day, are constructed from the prices of
Eurodollar futures contracts with decreasing expiration dates.  The
seller of a Eurodollar futures contract takes on an obligation to
deliver a deposit of one million US dollars to a bank account outside
the United States $i$ months from today. The price the buyer is
willing to pay for this contract depends on the prevailing interest
rate.  These contracts are traded on the Chicago Mercantile Exchange
and provide a way to lock in an interest rate. Eurodollar futures are
a liquid asset and are responsive to the Federal Reserve policy,
inflation, and economic indicators.

The data we study consist of 114 points per day; point $i$ corresponds
to the price of a contract with closing date $i$ months from today.
We consider 4 samples, each consisting of 100 days of this data:

{\em Sample 1:} curves from September 7, 1999 to January 27, 2000.

{\em Sample 2:} curves from January 24, 1997 to June 17, 1997.

{\em Sample 3:} curves from December 4, 1995 to April 24, 1996.

{\em Sample 4:} curves from March 6, 2001 to July 26, 2001.

\begin{figure}[ht!]

\begin{capf}\label{f:sample_mean_comparison}
{\small Sample means of the Eurodollar curves:
Left: sample 1 solid; sample 2 dashed.
Right: sample 3 solid; sample 4 dashed.}
\end{capf}

\centering

\includegraphics[width=7cm]{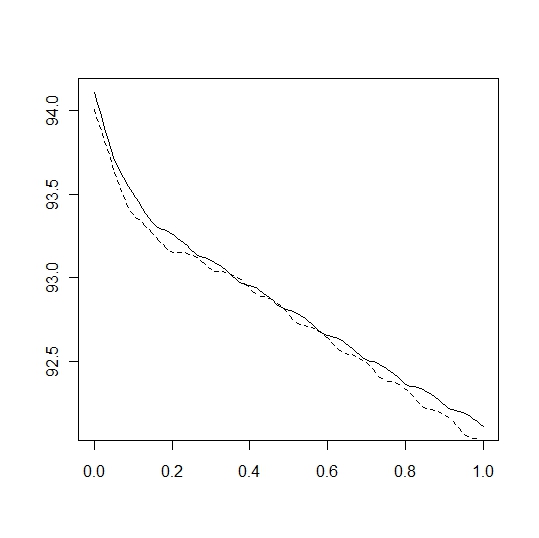}\hfill
\includegraphics[width=7cm]{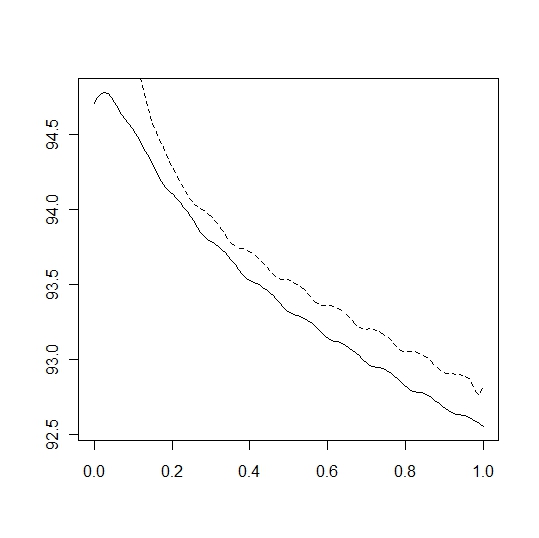}

\end{figure}

Figure~\ref{f:sample_mean_comparison} shows the sample mean functions for
the four samples. If a significance test does not reject $H_0$, we
can conclude that the expectations of the future
evolution of interest rates are the same for the two periods over
which the samples were taken. A rejection means that these expectations
are significantly different. As the analysis below reveals, we can
conclude that expectations of future interest rates were different in
Spring 1996 than in Summer 2001.

Table~\ref{tb:eurodollar} shows the P--values as a function of $p$
when the test is applied to samples 1 and 2,  and also when it is applied
to samples 3 and 4.  In both samples 1 and 2, $p=1$ explains more than
94\% of the variance, in both samples 3 and 4, $p=1$ explains
more than 84\% of the variance. Thus, following the recommendation
of Section \ref{s:te}, we use the P--values obtained with $p=1$.
They lead to the acceptance  of the null hypothesis of the equality
of mean functions for periods corresponding to samples 1 and 2,  and to
its rejection for periods corresponding to samples 3 and 4
(notice that the 0.81 in the bottom panel of Table~\ref{tb:eurodollar}
is 0.81\%).  These conclusions agree with a visual evaluation
of the sample mean functions in  Figure~\ref{f:sample_mean_comparison}.
They also confirm the observation made in Section \ref{ss:sim} that the
tests have very good power, as the curves in the right panel of
Figure~\ref{f:sample_mean_comparison} are not far apart. Both graphs in
Figure~\ref{f:sample_mean_comparison} give us an idea what kind of
differences in the sample mean functions are statistically significant,
and which are not.

\begin{table}

\begin{capt} \label{tb:eurodollar}

{\small P-values (in percent) of the statistics
applied to Eurodollar data;
samples 1 and 2 (top),  samples 3 and 4 (bottom).}
\end{capt}

\begin{center} \begin{tabular}{|c|rrrrr|} \hline
$p$      & 1  & 2 & 3 & 4 & 5 \\
\hline
$U^{(1)}$& 38.49 & 39.17 & 37.12 & 37.15 & 35.50\\
$U^{(2)}$& 38.49 & 68.52 & 0.00 & 0.00 & 0.00\\ \hline \end{tabular}\\

\vspace{2mm}

\begin{tabular}{|c|rrrrr|}
\hline
$p$      & 1  & 2 & 3 & 4 & 5\\
\hline
$U^{(1)}$& 0.81 & 0.23 & 0.10 & 0.01 & 0.07\\
$U^{(2)}$& 0.81 & 0.00 & 0.00 & 0.00 & 0.00\\ \hline \end{tabular}
\end{center} \end{table}

\subsection{Conclusions}
The simulations and data examples of this section show the tests
we propose enjoy good finite sample properties. Tests of this type
allow us to quantify statistical significance of conjectures made on the
basis of exploratory analysis. For example, the sample mean curves
in Figure~\ref{f:mean-flies} and the right panel of
Figure~\ref{f:sample_mean_comparison} look a bit different, but
a significance test allows us to state with confidence that
they correspond to different population mean functions.

In many procedures of functional data analysis, both exploratory and
inferential, the issue of choosing an optimal dimension reduction
parameter, like the $p$ in our setting, is delicate. Therefore,
procedures less sensitive to such a choice are preferable. From this
angle, the Monte Carlo test based $U^{(1)}$ is
preferable, as an inspection of Tables \ref{tb:mean-flies} and
\ref{tb:eurodollar} reveals. The test based on $U^{(2)}$ is
however easier to apply, and in our examples and simulations
leads to the same conclusions if $p$ is chosen according to the
cumulative variance rule.

The data examples of this section also show that the optimal value
of $p$ is typically a small single digit number, 1 or 2 in our examples.
Therefore, developing asymptotics as $p$ tends to infinity is not necessary,
and may, in fact, be misleading because for larger values of $p$ the tests
may yield counterintuitive results.

\section{Proofs of the results of
Section \ref{s:na}} \label{s:p-na}

{\sc Proof of Theorem \ref{t:2.1c}.}
The proof is done in two steps. First we show that $N^{-1/2}
\sum_{i=1}^N \varepsilon_i(t)$ is close to  $N^{-1/2}
\sum_{i=1}^N \varepsilon_{i,m}(t)$, if $m$ is sufficiently  large.
Then we establish (\ref{co-1}) for $m$-dependent functions for any $m\geq 1$.

As the first step, we show that
\beq \label{N1-1}
\limsup_{m\rightarrow \infty}\limsup_{N\rightarrow \infty}E\int  \lb N^{-1/2}
\sum_{i=1}^N\lp \varepsilon_i(t)-\varepsilon_{i,m}(t)\rp\rb^2dt=0,
\eeq
where the variables $\varepsilon_{i,m}$ are defined in (\ref{e:1.8a}).
By stationarity,
\begin{align*}
E&\lb \sum_{1\leq i \leq N}\lp \varepsilon_i(t)-\varepsilon_{i,m}(t)\rp\rb^2\\
&=\sum_{1\leq i \leq N}\sum_{1\leq j \leq N} E\lp \varepsilon_i(t)-\varepsilon_{i,m}(t)\rp\lp \varepsilon_j(t)-\varepsilon_{j,m}(t)\rp\\
&=NE\lp \varepsilon_0(t)-\varepsilon_{0,m}(t)\rp^2+2\sum_{1\leq i<j\leq N}
E\lp \varepsilon_i(t)-\varepsilon_{i,m}(t)\rp\lp \varepsilon_j(t)-\varepsilon_{j,m}(t)\rp.
\end{align*}
In the proof, we will repeatedly use independence relations which
follow from representations (\ref{e:1.6}) and (\ref{e:1.8a}). First
observe that  if $j>i$, then
$(\varepsilon_i, \varepsilon_{i,m})$ is independent of
$\varepsilon_{j,j-i}$ because
\bd
\eg_{j, j-i} = f(\dg_j, \ldots, \dg_{i+1}, \dg_{j, i}^{(j-i)},
\dg_{j, i-1}^{(j-i)}, \ldots).
\ed
Consequently,
$E\lp\varepsilon_i(t)-\varepsilon_{i,m}(t)\rp \varepsilon_{j,j-i}(t)=0$,
and so
\bd
\sum_{1\leq i<j\leq N}
E\lp \varepsilon_i(t)-\varepsilon_{i,m}(t)\rp \varepsilon_j(t)
=\sum_{1\leq i<j\leq N}
E\lp \varepsilon_i(t)-\varepsilon_{i,m}(t)\rp\lp
\varepsilon_j(t)-\varepsilon_{j,j-i}(t)\rp.
\ed
Using the Cauchy-Schwarz inequality and (\ref{e:1.9}), we conclude
\begin{align*}
\biggl| \int \sum_{1\leq i<j\leq N}&
E\lp \varepsilon_i(t)-\varepsilon_{i,m}(t)\rp\lp \varepsilon_j(t)-\varepsilon_{j,j-i}(t)\rp dt\biggl|\\
&\leq \sum_{1\leq i<j\leq N}\int \lb E\lp \varepsilon_i(t)-\varepsilon_{i,m}(t)\rp^2\rb^{1/2}\lb E\lp \varepsilon_j(t)-\varepsilon_{j,j-i}(t)\rp^2\rb^{1/2}dt\\
&\leq  \sum_{1\leq i<j\leq N}
\lb \int E\lp \varepsilon_i(t)-\varepsilon_{i,m}(t)\rp^2dt\rb^{1/2}
\lb \int E\lp \varepsilon_j(t)-\varepsilon_{j,j-i}(t)\rp^2
dt\rb^{1/2}\\
&= \sum_{1\leq i<j\leq N}
\lb \int E\lp \varepsilon_0(t)-\varepsilon_{0,m}(t)\rp^2dt\rb^{1/2}
\lb \int E\lp \varepsilon_0(t)-\varepsilon_{0,j-i}(t)\rp^2
dt\rb^{1/2}\\
&\leq   N\lb \int
E\lp \varepsilon_{0,m}(t)- \varepsilon_{0}(t)\rp^2dt\rb^{1/2}
\sum_{k\geq 1}\lb
\int\lp \varepsilon_0(t)-\varepsilon_{0,k}(t)\rp^2dt\rb^{1/2}.
\end{align*}
Hence
$$
\limsup_{m\to \infty}\limsup_{N\to \infty}\frac{1}{N}\left|
\int \sum_{1\leq i<j\leq N}
E\lb\lp \varepsilon_{i,m}(t)-\varepsilon_{i}(t)\rp \varepsilon_j(t)\rb dt
\right|=0.
$$
Similar arguments give
$$
\limsup_{m\to \infty}\limsup_{N\to \infty}\frac{1}{N}\left|\int \sum_{1\leq i<j\leq N}E\lb\lp \varepsilon_{i,m}(t)-\varepsilon_{i}(t)\rp \varepsilon_{j,m}(t)\rb dt
\right|=0.
$$
Completing the verification of (\ref{N1-1}).

The next the step is to show that $N^{-1/2} \sum_{1 \le i \le N} \eg_{i,m}$
converges to a Gaussian process $Z_m$ with covariances defined analogously
to (\ref{e:cts}).   Recall that for every integer
$m\geq 1$,  $\{\varepsilon_{i,m}\}$ is an
$m$--dependent sequence of functions. To lighten the notation,
in the remainder of the proof, we fix $m$ and denote sequence
$\{ \eg_{i,m} \}$ by $\{ \eg_{i} \}$, so $\{ \eg_{i} \}$ is
now $m$--dependent.

Let $K>1$ be an integer and
$\psi_i$ be an orthonormal basis determined by the eigenfunctions of
$E\varepsilon(t)\varepsilon(s)$. The corresponding eigenvalues are
denoted by $\nu_i$. Then, by the Karhunen-Lo\'eve expansion, we have
$$
\varepsilon_i(t)
=\sum_{\ell\geq 1}\lip \varepsilon_i, \psi_\ell \rip \psi_\ell(t).
$$
Next we define
$$
\varepsilon_i^{(K)}(t)=\sum_{1\leq \ell \leq K}\lip \varepsilon_i, \psi_\ell \rip \psi_\ell(t).
$$
By the triangle inequality we have that
\begin{align*}
&\lbr E\int\lb\sum_{1\leq i \leq N}\lp(\varepsilon_i(t)-\varepsilon_i^{(K)}(t)\rp\rb^2dt\rbr^{1/2}\\
&\hspace{.5 cm}\leq \lbr E\int\lb\sum_{i\in V(0)}\lp(\varepsilon_i(t)-\varepsilon_i^{(K)}(t)\rp\rb^2dt\rbr^{1/2}+\ldots\\
&\hspace{1.5 cm} +\lbr E\int\lb\sum_{ i\in V(m-1)}\lp(\varepsilon_i(t)-\varepsilon_i^{(K)}(t)\rp\rb^2dt\rbr^{1/2},
\end{align*}
where $V(k)=\{i:\; 1\leq i \leq N, i=k\;(\mbox{mod}\; m)\}, 0\leq k \leq m-1$. Due to the $m$ dependence of the sequence $\{\varepsilon_i\}$, $\sum_{ i\in V(k)}(\varepsilon_i(t)-\varepsilon_i^{(K)}(t))$ is a sum of independent, identically distributed random variables, and thus we get that
$$
E\int\lb\sum_{ i\in V(m-1)}\lp(\varepsilon_i(t)-\varepsilon_i^{(K)}(t)\rp\rb^2dt\leq
N\sum_{\ell\geq K}E\lip X_0, \psi_\ell \rip^2.
$$
Utilizing
$$
\lim_{K\to \infty}\sum_{\ell\geq K}E\lip X_0, \psi_\ell \rip^2=0
$$
we conclude that for any $x>0$
$$
\limsup_{K\to \infty}\limsup_{N\to \infty}P\lbr\int\lb\frac{1}{N^{1/2}}\sum_{1\leq i \leq N}\lp(\varepsilon_i(t)-\varepsilon_i^{(K)}(t)\rp\rb^2dt>x \rbr=0.
$$
The sum of the $\varepsilon_i^{(K)}$'s can be written as
$$
\frac{1}{N^{1/2}}\sum_{1\leq i \leq N}\varepsilon_i^{(K)}(t)=\sum_{1\leq \ell \leq K}\psi_\ell(t)\frac{1}{N^{1/2}}\sum_{1\leq i \leq N}\lip \varepsilon_i, \psi_{\ell}\rip.
$$
Next, we use the central limit theorem for stationary $m$-dependent sequences
of random vectors
(see \citetext{lehmann:1999}
 and the Cram\'er-Wold theorems in \citetext{dasgupta:2008},
pages 9 and 120)) and  and get that
$$
\lbr \frac{1}{N^{1/2}}\sum_{1\leq i \leq N}\lip \varepsilon_i, \psi_{\ell}\rip, 1\leq \ell \leq K\rbr^T \convd {\bf N}_K({\bf 0}, \bDg_K),
$$
where ${\bf N}_K({\bf 0}, \bDg_K)$ is a $K$-dimensional
normal random variable with zero mean and covariance matrix
$
\bDg_K=\mbox{diag}(\nu_1,\ldots,\nu_K).
$
Thus we proved that for all $K>1$
$$
N^{-1/2}\sum_{1\leq i \leq N}\varepsilon_i^{(K)}(t)
\convd \sum_{1\leq \ell \leq K}\nu_\ell^{1/2}N_\ell\psi_\ell(t)
\;\;\;\mbox{in}\;\;L^2,
$$
where $N_i, i\geq 1$ are independent standard normal random variables.
It is easy to see that
$$
\int
\lp\sum_{ K<\ell <\infty}\nu_\ell^{1/2}N_\ell\psi_\ell(t)\rp^2dt=
\sum_{ K<\ell <\infty}\nu_\ell N_\ell^2
\convP 0,
$$
as $K\to \infty$. Thus we have the convergence of
$N^{-1/2}\sum_{1\leq i \leq N}
\varepsilon_i$ for any $m$ and therefore the proof of the Theorem
is now complete.

\vspace{4mm}

{\sc Proof of Theorem \ref{t:3.1c}.}
\label{s:PD}
First we reduce (\ref{e:6.1}) to (\ref{e:6.4}).
Then,  we reduce (\ref{e:6.4}) to (\ref{e:6.9}).

Since
\bd
\hat\ga_0(t,s) = \frac{1}{N} \sum_{i=1}^N
 \lp X_i(t) - \mu(t) \rp \lp X_i(s) - \mu(s) \rp
- \lp \bar X_N(t) - \mu(t) \rp  \lp \bar X_N(s) - \mu(s) \rp,
\ed
we obtain that
\begin{align*}
\dint &\lbr \hat\ga_0(t,s) - E[\eg_0(t)\eg_0(s)] \rbr^2dtds\\
&\le 4 \dint \lbr \frac{1}{N} \sum_{i=1}^N
 \lp X_i(t) - \mu(t) \rp \lp X_i(s) - \mu(s) \rp
- E[\eg_0(t)\eg_0(s)]
 \rbr^2 dtds\\
&\hspace{1 cm} + 4 \lp \int \lp \bar X_N(t) - \mu(t) \rp^2dt \rp^2\\
&=o_P(1),
\end{align*}
using the ergodic theorem for random variables in a Hilbert space.

Next observe that
\begin{align*}
\hat\ga_i(t,s)
&= \frac{1}{N} \sum_{j=i+1}^N \eg_j(t) \eg_{j-i}(s) \\
&\ \ + \frac{N-i}{N}\bar \eg_N(t) \bar \eg_N(s)
- \lp \frac{1}{N} \sum_{j=i+1}^N \eg_j(t)\rp \bar \eg_{N}(s)
- \bar \eg_{N}(t) \lp \frac{1}{N} \sum_{j=i+1}^N \eg_{j-i}(s)\rp.
\end{align*}

Therefore, setting,
\bd
 \bar \ga_i(t,s) = \frac{1}{N} \sum_{j=i+1}^N\eg_j(t) \eg_{j-i}(s).
\ed
we obtain
\begin{align*}
\sum_{i=1}^{N-1} K\lp\frac{i}{h}\rp \hat\ga_i(t,s)
&= \sum_{i=1}^{N-1} K\lp\frac{i}{h}\rp \bar \ga_i(t,s)\\
&\ \ - \sum_{i=1}^{N-1} K\lp\frac{i}{h}\rp
\lbr \frac{1}{N} \sum_{j=i+1}^N \eg_j(t)\rbr \bar\eg_{N}(s)\\
&\ \ - \sum_{i=1}^{N-1} K\lp\frac{i}{h}\rp
\lbr \frac{1}{N} \sum_{j=i+1}^N \eg_{j-i}(s)\rbr \bar\eg_{N}(t)\\
&\ \ + \sum_{i=1}^{N-1} K\lp\frac{i}{h}\rp
\frac{N-i}{N}\bar\eg_{N}(t) \bar \eg_{N}(s).
\end{align*}

gud

By stationarity we conclude that for any $1\le i \le N$,
\begin{align*}
E \int \lp \frac{1}{N} \sum_{j=i+1}^N \eg_j(t)\rp^2 dt
&= \frac{1}{N^2} \sum_{j=i+1}^N  \int E \eg_0^2(t)dt
+ \frac{2}{N^2} \sum_{j=i+1}^N (N-j)\int E\eg_0(t) \eg_j(t)dt\\
&\leq \frac{1}{N} \int E \eg_0^2(t)dt
+ \frac{4}{N} \sum_{j=1}^\infty \left | \int E\eg_0(t) \eg_j(t)dt \right |.
\end{align*}
Since $\eg_0$ and $\eg_{j,j}$ are independent, we get by
(\ref{e:1.9})
\begin{align*}
\sum_{j=1}^\infty \left | \int E\eg_0(t) \eg_j(t)dt \right |
&= \sum_{j=1}^\infty
\left | \int E\eg_0(t) \lp \eg_j(t) - \eg_{j,j}(t) \rp dt \right | \\
&\le \sum_{j=1}^\infty
E\lbr
\lp \int \eg_0^2(t) dt \rp^{1/2}
\lp \int \lp \eg_j(t) - \eg_{j,j}(t) \rp^2 dt  \rp^{1/2}
\rbr\\
&\le \lp E \int \eg_0^2(t) dt \rp^{1/2}
\sum_{j=1}^\infty
\lp \int E \lp \eg_j(t) - \eg_{j,j}(t) \rp^2 dt  \rp^{1/2} \\
&< \infty.
\end{align*}
Thus, we have
\bd
\max_{1 \le i \le N} E \int \lp \frac{1}{N} \sum_{j=i+1}^N \eg_j(t) \rp^2 dt
= O(1).
\ed
Consequently, using the triangle inequality,
\begin{align*}
&E
\lbr
\dint
\lp
\sum_{i=1}^{N-1} K \lp \frac{i}{h} \rp
\lb \frac{1}{N} \sum_{j=i+1}^N \eg_j(t) \rb \bar\eg_N(s)
\rp^2 dtds
\rbr^{1/2} \\
&\le \sum_{i=1}^{N-1} \lmo K\lp \frac{i}{h} \rp  \rmo
E
\lbr
\lp \int \lb \frac{1}{N} \sum_{j=i+1}^N \eg_j(t) \rb^2 dt \rp^{1/2}
\lp \int \bar\eg_N^2(s) ds \rp^{1/2}
\rbr \\
&\le \sum_{i=1}^{N-1} \lmo K\lp \frac{i}{h} \rp  \rmo
\lp E \int \lb \frac{1}{N} \sum_{j=i+1}^N \eg_j(t) \rb^2 \rp^{1/2}
\lp \int \bar\eg_N^2(s) ds \rp^{1/2}\\
&= \frac{h}{N} O(1) = o(1),
\end{align*}
on account of (\ref{e:3.5}).

Hence to establish (\ref{e:6.1}) , it is enough to prove that
\beq \label{e:6.4}
\dint \lp  \sum_{i=1}^{N-1} K\lp\frac{i}{h}\rp
\bar \ga_i(t,s) - c_1(t,s) \rp^2dtds = o_P(1),
\eeq
where
\bd
c_1(t,s) = \sum_{i\ge 1} E [\eg_0(s) \eg_i(t)].
\ed

Let $\lbr \eg_{n,m}, \ -\infty < n < \infty \rbr$ be the
random variables defined in (\ref{e:1.8a}), where $m$ is a fixed
number. Let
\bd
\tilde\ga_{i,m}(t,s) = \frac{1}{N} \sum_{j=i+1}^{N}
\eg_{j,m}(t) \eg_{j-i,m}(s).
\ed
We show that for every $m\ge 1$,
\beq \label{e:6.6}
\dint \lp  \sum_{i=1}^{N-1} K\lp\frac{i}{h}\rp
\tilde\ga_{i,m}(t,s) - c_1^{(m)}(t,s) \rp^2dtds = o_P(1),
\eeq
where
\bd
c_1^{(m)}(t,s) = \sum_{i=1}^\infty E [ \eg_{1,m}(s)\eg_{i+1,m}(t)].
\ed
We also note that (\ref{e:1.8a}) and (\ref{e:1.9}) imply
\beq \label{e:6.7}
\lim_{m\to\infty}
\dint \lp c_1^{(m)}(t,s) - c_1(t,s) \rp^2dtds = 0.
\eeq
Since $\lbr \eg_{n,m}, \ -\infty < n < \infty \rbr$ is an $m$--dependent
sequence,
\bd
c_1^{(m)}(t,s) = \sum_{i=1}^m E [ \eg_{1,m}(s)\eg_{i+1,m}(t)].
\ed
Using  (\ref{e:3.1}), (\ref{e:3.2}) and  (\ref{e:3.5}), we get
\bd
\max_{1 \le i \le m} \lmo  K\lp\frac{i}{h}\rp -1 \rmo \to 0, \ \ \
{\rm as}\ N\to \infty.
\ed
By the ergodic theorem,
\bd
\dint \lp \tilde\ga_{i,m}(t,s) - E[ \eg_{1,m}(s)\eg_{i+1,m}(t)]\rp^2dtds
=o_P(1),
\ed
for any fixed $i$.
Hence (\ref{e:6.6}) is proven, once we have shown that
\beq \label{e:6.8}
\dint \lp  \sum_{i=m+1}^{N-1} K\lp\frac{i}{h}\rp\tilde\ga_{i,m}(t,s)
\rp^2dtds=o_P(1).
\eeq

It is easy to see that
\begin{align*}
E &\dint \lp  \sum_{i=m+1}^{N-1} K\lp\frac{i}{h}\rp\tilde\ga_{i,m}(t,s)
\rp^2dtds
\\
&= \dint
\lp \frac{1}{N^2}
\sum_{i=m+1}^h\sum_{\ell=m+1}^h\sum_{k=i+1}^{N-1}\sum_{n=\ell+1}^{N-1}
K\lp\frac{i}{h}\rp K\lp\frac{\ell}{h}\rp
E\lb \eg_{k,m} \eg_{k-i,m}\eg_{n,m}\eg_{n-\ell,m} \rb
\rp,
\end{align*}
provided $h\le N-1$. The sequence  $\lbr \eg_{n,m}, \ -\infty < n < \infty \rbr$ is an $m$--dependent, and therefore
 $\eg_{k,m}$ and $\eg_{k-i,m}$ are independent, since $i\ge m+1$.
Similarly, $\eg_{n,m}$ and $\eg_{n-\ell,m}$ are independent.
Hence the number of terms when
$E\lb \eg_{k,m} \eg_{k-i,m}\eg_{n,m}\eg_{n-\ell,m} \rb$ is
different from zero is $O(Nh)$. Consequently,
\bd
E \dint \lp  \sum_{i=m+1}^{N-1} K\lp\frac{i}{h}\rp\tilde\ga_{i,m}(t,s)
\rp^2dtds = O\lp \frac{h}{N}\rp= o(1).
\ed
This completes the verification of( \ref{e:6.8}).

Next we show that for all $\epsilon > 0$,
\begin{align} \label{e:6.9}
\lim_{m\to\infty}& \limsup_{N\to \infty}
P\lbr
\dint
\lp
 \sum_{i=1}^{N-1} K\lp\frac{i}{h}\rp
\lb
\bar\ga_{i}(t,s) - \tilde\ga_{i,m}(t,s)
\rb
\rp^2dtds > \epsilon
\rbr = 0.
\end{align}
Using the definitions of the covariances $\bar\ga_{i}(t,s)$ and
$\tilde\ga_{i,m}(t,s)$, we consider the decompositions
\begin{align*}
\frac{1}{N}  \sum_{i=1}^{N-1}& K\lp\frac{i}{h}\rp
\sum_{j=i+1}^N
\lb
\eg_j(t) \eg_{j-i}(s) - \eg_{j,m}(t) \eg_{j-i,m}(s)
\rb\\
&= \frac{1}{N} \lbr \sum_{i=1}^m +   \sum_{i=m+1}^h
\rbr K\lp\frac{i}{h}\rp \sum_{j=i+1}^N
\lb
\eg_j(t) \eg_{j-i}(s) - \eg_{j,m}(t) \eg_{j-i,m}(s)
\rb
\end{align*}
and
\begin{align*}
\eg_j(t) &\eg_{j-i}(s) - \eg_{j,m}(t) \eg_{j-i,m}(s)
= \lp \eg_j(t) - \eg_{j,m}(t)\rp\eg_{j-i}(s)
+ \lp \eg_{j-i}(s) - \eg_{j-i,m}(s) \rp \eg_{j,m}(t).
\end{align*}
Clearly,
\begin{align*}
 &\hspace{-.5 cm}\lbr
\dint
\lp
\frac{1}{N}  \sum_{i=1}^m K\lp\frac{i}{h}\rp \sum_{j=i+1}^N
\lp \eg_j(t) - \eg_{j,m}(t)\rp\eg_{j-i}(s)
\rp^2dtds
\rbr^{1/2}
\\
&\le \frac{1}{N}\sum_{i=1}^m \lmo K\lp\frac{i}{h}\rp\rmo
\lbr \int \lp \eg_j(t) - \eg_{j,m}(t)\rp^2 dt\rbr^{1/2}
\lbr \int \eg_{j-i}^2(s)ds\rbr^{1/2},
\end{align*}
so, by (\ref{m1}),
\begin{align*}
&\hspace{-.5 cm}E \lbr
\dint
\lp
\frac{1}{N}  \sum_{i=1}^m K\lp\frac{i}{h}\rp \sum_{j=i+1}^N
\lp \eg_j(t) - \eg_{j,m}(t)\rp\eg_{j-i}(s)
\rp^2dtds
\rbr^{1/2}\\
&\le m
\lbr
E \int \lp \eg_0(t) - \eg_{0,m}(t)\rp^2 dt E \int \eg_{0}^2(s)ds
\rbr^{1/2} \\
&\le A m \lbr
E \int \lp \eg_0(t) - \eg_{0,m}(t)\rp^2 dt\rbr^{1/2} \to 0,
\ \ \ {\rm as}\ \  m \to \infty,
\end{align*}
according to (\ref{m1}), where $A$ is a constant.

Next we use the decomposition
\bd
\eg_j(t) \eg_{j-i}(s)
= \eg_{j,i}(t) \eg_{j-i}(s) + \lp \eg_j(t)
- \eg_{j,i}(t)\rp \eg_{j-i}(s)
\ed
to get
\begin{align*}&\hspace{-.5 cm}
\lbr
\dint
\lp
\frac{1}{N}  \sum_{i=m+1}^h K\lp\frac{i}{h}\rp \sum_{j=i+1}^N
\lp \eg_j(t) - \eg_{j,i}(t)\rp\eg_{j-i}(s)
\rp^2dtds
\rbr^{1/2}
\\
&\le \frac{1}{N}\sum_{i=m+1}^\infty  \sum_{j=i+1}^N
\lbr \int \lp \eg_j(t) - \eg_{j,i}(t)\rp^2 dt\rbr^{1/2}
\lbr \int \eg_{j-i}^2(s)ds\rbr^{1/2}.
\end{align*}
Therefore, by (\ref{e:1.8}) and (\ref{e:1.9}), we have
\begin{align*}
&E \lbr
\dint
\lp
\frac{1}{N}  \sum_{i=m+1}^h K\lp\frac{i}{h}\rp \sum_{j=i+1}^N
\lp \eg_j(t) - \eg_{j,i}(t)\rp\eg_{j-i}(s)
\rp^2dtds
\rbr^{1/2}\\
&
\le \frac{1}{N}\sum_{i=m+1}^\infty \sum_{j=i+1}^N
E \lb \lbr \int \lp \eg_j(t) - \eg_{j,i}(t)\rp^2 dt\rbr^{1/2}
\lbr \int \eg_{j-i}^2(s)ds\rbr^{1/2} \rb\\
&
\le \frac{1}{N}\sum_{i=m+1}^\infty \sum_{j=i+1}^N
\lb  \int \lp \eg_j(t) - \eg_{j,i}(t)\rp^2 dt \rb^{1/2}
\lb \int \eg_{j-i}^2(s)ds \rb^{1/2}\\
&
\le A \sum_{i=m+1}^\infty
\lb  \int \lp \eg_0(t) - \eg_{0,i}(t)\rp^2 dt \rb^{1/2}
\to 0, \ \ \ {\rm as} \ \ m\to \infty.
\end{align*}

We have  shown so far that for any $\epsilon > 0$,
\bd
\lim_{m\to\infty} \limsup_{N\to \infty}
P \lbr \dint
\lp
\frac{1}{N} \sum_{i=m+1}^h K\lp\frac{i}{h}\rp \sum_{j=i+1}^N
\eg_{j}(t) \eg_{j-i}(s)
\rp^2 dtds > \epsilon \rbr = 0.
\ed
Similar arguments give
\bd
\lim_{m\to\infty} \limsup_{N\to \infty}
P\lbr \dint \lp
\frac{1}{N} \sum_{i=m+1}^h K\lp\frac{i}{h}\rp \sum_{j=i+1}^N
\eg_{j, m}(t) \eg_{j-i,m}(s)
\rp^2 dtds > \epsilon \rbr = 0.
\ed

This completes the verification of (\ref{e:6.9}), so (\ref{e:6.4}) is
proven.

\section{Proofs of the results
of Section \ref{s:te}} \label{s:p-te}

In the proofs, we will often use the relations
\beq \label{e:1.13}
\max_{1 \le i \le p } |\hat\la_i - \la_i|\convP 0, \ \
{\rm and} \ \
\max_{1 \le i \le p } \lnorm \hat\fg_i - \hc_i \fg_i\rnorm \convP 0
\ \ \
{\rm as} \ \min(M,N) \to \infty,
\eeq
where
$
\hc_i = {\rm sign} \lp \lip \hat\fg_i, \fg_i \rip \rp,
$
Analogous  relations  have been extensively used for the eigenvalues and
the eigenfunctions of the empirical and population covariance operators,
see \citetext{bosq:2000}, \citetext{gabrys:kokoszka:2007},
 \citetext{horvath:huskova:kokoszka:2009},
\citetext{panaretos:2010}, among many others, but
they hold in much greater generality, see Chapter 2 of \citetext{HKbook}.
Under (\ref{e:1.12}), they hold for the eigenelements of the
operators $\hat D$ and $D$ defined in Section~\ref{s:te}.

\begin{proof}[{\sc  Proof of Theorem \ref{t:1.1}:}] By Theorem~\ref{t:2.1c},
assumptions (\ref{e:1.6})--(\ref{e:1.9}) imply that, as $N\to \infty$
and $M\to \infty$
\bd
\lp N^{-1/2} \sum_{1\le i \le N}\eg_i,\  M^{-1/2} \sum_{1\le j \le N}\eg_j^*\rp
\convd \lp \Gg^{(1)},\Gg^{(2)} \rp,
\ed
where $\Gg^{(1)}$ and $\Gg^{(2)}$ are independent $L^2$--valued
mean zero Gaussian processes with covariances
\bd
E \lb \Gg^{(1)}(t)\Gg^{(1)}(s) \rb = c(t,s), \ \ \
 E \lb \Gg^{(2)}(t)\Gg^{(2)}(s) \rb = c^*(t,s).
\ed

Hence
\beq \label{e:4.1}
\frac{NM}{N+M} \int
\lp  N^{-1} \sum_{1\le i \le N}\eg_i(t)
-  M^{-1} \sum_{1\le j \le M}\eg_j^*(t)\rp^2 dt
\convd \int \Gg(t) dt,
\eeq
where
$
\Gg(t) = (1-\theta)^{1/2} \Gg^{(1)}(t) + \theta^{1/2} \Gg^{(2)}(t).
$
The conclusion of Theorem \ref{t:1.1} now follows.
\end{proof}

\vspace{2mm}

\begin{proof}[{\sc Proof of Theorem \ref{t:1.2}:}] By the ergodic theorem
in a Hilbert space, $|| \bar X_N -\mu || = o_P(1)$
and $|| \bar X_M -\mu^* || = o_P(1)$, which imply the result.
\end{proof}

\vspace{2mm}

\begin{proof}[{\sc Proof of Theorem \ref{t:2.1}:}]
 Under assumptions (\ref{e:1.8}) and (\ref{e:1.9}), we have, jointly,
\bd
\lbr
N^{-1/2} \sum_{1\le i \le N} \lip \eg_i, \fg_k \rip, \ 1 \le k \le p
\rbr
\convd {\bf N}_p^{(1)}\lp \bzero, \bQ^{(1)} \rp
\ed
and
\bd
\lbr
M^{-1/2} \sum_{1\le j \le M} \lip \eg_j^*, \fg_k \rip, \ 1 \le k \le p
\rbr
\convd {\bf N}_p^{(2)}\lp \bzero, \bQ^{(2)} \rp,
\ed
where ${\bf N}_p^{(1)}\lp \bzero, \bQ^{(1)}\rp$ and
${\bf N}_p^{(2)}\lp \bzero, \bQ^{(2)}\rp$  are independent
$p$--dimensional normal random vectors. Since
\bd
a_k =   N^{-1} \sum_{1\le i \le N} \lip \eg_i, \fg_k \rip
-  M^{-1} \sum_{1\le j \le M} \lip \eg_j^*,  \fg_k \rip. \ \ \
1 \le i \le p,
\ed
(\ref{e:2.2}) holds with
$\bQ = (1-\theta)\bQ^{(1)} + \theta  \bQ^{(2)}$.

Next we observe that the matrix
$\bQ = \lp Q(k,\ell), \ 1 \le k,\ell \le p\rp$ satisfies
\bd
 Q(k,\ell) =\dint d(t,s) \fg_k(t) \fg_\ell(s) dtds
= \la_k \dg_{k\ell},
\ed
where $\dg_{ij}$ is Kronecker's delta, using the fact that
the $\fg_i$ are orthonormal eigenfunctions and the $\la_i$ are
the corresponding eigenvalues.

Recall now relations (\ref{e:1.13}), and observe that
\begin{align*}
\hat a_k
&= \hc_k
\lp  N^{-1} \sum_{1\le i \le N} \lip \eg_i, \fg_k \rip
-  M^{-1} \sum_{1\le j \le M} \lip \eg_j^*,  \fg_k \rip
\rp \\
& \hspace{1.5 cm} +
\lip
N^{-1} \sum_{1\le i \le N} \eg_i
- M^{-1} \sum_{1\le j \le M}  \eg_j^*\ , \
\hat\fg_k - \hc_k \fg_k
\rip.
\end{align*}
By (\ref{e:1.13}) and (\ref{e:4.1}),
\begin{align*}
\lp \frac{NM}{N+M}\rp^{1/2}&
\lmo \lip
\frac{1}{N} \sum_{1\le i \le N} \eg_i
- \frac{1}{M} \sum_{1\le j \le M}  \eg_j^*\ , \
\hat\fg_k - \hc_k \fg_k
\rip \rmo\\
&\le  U_{N,M}^{1/2}
\lp \int (\hat\fg_k - \hc_k \fg_k)^2 dt  \rp^{1/2} \\
&= o_P(1).
\end{align*}
Hence the result follows immediately from (\ref{e:2.2}) and the
diagonality of $\bQ$.
\end{proof}

\vspace{2mm}

\begin{proof}[{\sc Proof of Theorem \ref{t:2.2}:}] By the ergodic
theorem $a_i \convP \lip \mu - \mu^*, \fg_i \rip, \ 1 \le i \le p$.
Since relations (\ref{e:1.13}) hold also under $H_A$, the result is proven.
\end{proof}

\vspace{2mm}

\noindent{\bf Acknowledgements.} The research was partially supported by
NSF grants DMS-0905400 and DMS-0931948.

\bibliographystyle{oxford3}
\renewcommand{\baselinestretch}{0.9}
\small

\bibliography{Z:/TEX/kokoszka}

\end{document}